\documentclass[preprint,12pt,1p]{elsarticle}
\makeatletter
\def\ps@pprintTitle{%
  \let\@oddhead\@empty
  \let\@evenhead\@empty
  \let\@oddfoot\@empty
  \let\@evenfoot\@oddfoot
}
\makeatother
\usepackage{graphics}
\usepackage{epsfig}
\usepackage{mathptmx}
\usepackage{amsmath}
\usepackage{amssymb}
\usepackage{hyperref}
\usepackage{fullpage}
\usepackage{amsthm}
\usepackage{geometry}
\theoremstyle{plain}

\theoremstyle{definition}

\theoremstyle{example}

\theoremstyle{remark}

\newlength{\defbaselineskip}
\journal{Computer Aided Geometric Design}
\numberwithin{equation}{section}
\begin{document}
\begin{frontmatter}
\title{Laplace transforms based some novel integrals via  hypergeometric technique}
\author{M. I. Qureshi}
\author{*Showkat Ahmad Dar}
\ead{showkat34@gmail.com}
\address{E-Mail: miqureshi\_delhi@yahoo.co.in and showkat34@gmail.com\\Department of Applied Sciences and Humanities , \\Faculty of  Engineering and Technology, \\Jamia Millia Islamia ( Central University), New Delhi, 110025, India.}
\cortext[cor1]{Corresponding author}
\begin{abstract}
  In this paper, we obtain the analytical solutions of Laplace transforms based some novel integrals with suitable convergence conditions, by using hypergeometric approach (some algebraic properties of Pochhammer symbol and  classical summation theorems of hypergeometric series ${}_{2}F_{1}(1)$, ${}_{2}F_{1}(-1)$ , ${}_{4}F_{3}(-1)$) . Also, we obtain the Laplace transforms of arbitrary powers of some finite series containing hyperbolic sine and cosine functions having different arguments, in terms of hypergeometric  and  Beta functions. Moreover, Laplace transforms of even and odd positive integral powers of sine and cosine functions with different arguments, and their combinations of the product (taking two, three, four functions at a time), are obtained. In addition, some special cases are yield from the main results.
 \\
 \\
\textit{2010 AMS Classification: 33C05; 33C20; 44A10; 33B15 }
 \end{abstract}
\end{frontmatter}
\section{Introduction and Preliminaries}
For the sake of conciseness of this paper, we use the following notations \\
$~~~~~~~~~~~\mathbb{N}:=\{1,2,...\};~~~~~~\mathbb{N}_{0}:=\mathbb{N}\cup\{0\};~~~~~~\mathbb{Z}_{0}^{-}:=\mathbb{Z}^{-}\cup\{0\}=\{0,-1,-2,-3,...\},$\\
where the symbols $\mathbb{N}$ and $\mathbb{Z}$ denote the set of natural numbers and integers; as usual, the symbols $\mathbb{R}$ and $\mathbb{C}$ denote the set of real and complex numbers.\\
$~~~~~$In the table of  Gradshteyn and Ryzhik  \cite[section 3.5]{G} on the definite integrals, some integrand contains quotient of the classical hyperbolic functions like as $\cosh(x)$ and $\sinh(x)$, which can be expressed in terms of the exponential functions, defined by
\begin{equation}\label{FAR1}
\cosh(x)=\frac{e^{x}+e^{-x}}{2}~~and~~\sinh(x)=\frac{e^{x}-e^{-x}}{2}\newline.
\end{equation}
The classical Beta function $B(\alpha,\beta)$ \cite[p. 26, eq.(48)]{s1}, is defined by
\begin{equation}\label{FAR2}
B(\alpha, \beta)=\begin{cases} \displaystyle{\int_{0}^{1}t^{\alpha-1}(1-t)^{\beta-1}dt}, & ~~~~\Re(\alpha)>0,~ \Re(\beta)>0,\\
\\
\displaystyle{\frac{\Gamma(\alpha)\Gamma(\beta)}{\Gamma(\alpha+\beta)}}, &  \Re(\alpha)<0,~ \Re(\beta)<0;~\quad \alpha,\beta\in\mathbb{C}\backslash \mathbb{Z}_{0}^{-}. \\
 \end{cases}
 \end{equation}
Gauss hypergeometric series ${}_{2}F_{1}(\cdot)$ \cite[p.29, Eq.(4)]{s1}, is defined by
 \begin{equation}\label{FAR3}
 {}_{2}F_{1}\left(\begin{array}{lll}a,~b~;\\~~~~ d~;\end{array} z \right)=
\sum_{n=0}^{\infty}\frac{(a)_{n}(b)_{n}}{(d)_{n}}\frac{z^{n}}{n!},
\end{equation}
~~~~~where $~|z|<1;~a,b\in\mathbb{C};~~d\in\mathbb{C}\backslash \mathbb{Z}_{0}^{-}.$\\
The notation $(\lambda)_{\upsilon}~(\lambda, \upsilon\in\mathbb{C})$ denotes the Pochhammer's symbol (or the shifted factorial,
 since $(1)_{n}=n! )$ is defined, in general, by
\begin{equation}\label{FAR5}
(\lambda)_{\upsilon}:=\frac{\Gamma(\lambda+\upsilon)}{\Gamma(\lambda)}=\begin{cases} 1 , \quad~~~~~~~~~~~~~ (\upsilon=0~;~\lambda\in\mathbb{C}\backslash\{0\}) \\ \lambda(\lambda+1)...(\lambda+n-1),  \quad (\upsilon=n\in\mathbb{N}~;~\lambda\in\mathbb{C}). \\
 \end{cases}
 \end{equation}
A natural generalization of Gauss hypergeometric series ${}_{2}F_{1}$ is the general hypergeometric series ${}_{p}F_{q}$  
 with $p$ numerator parameters $\alpha_{1},  ... , \alpha_{p}$ and $q$ denominator parameters $\beta_{1}, ..., \beta_{q}$. It is defined by\\
\begin{equation}\label{FAR4}
{}_{p}F_{q}\left(\ \begin{array}{lll}\alpha_{1},...,\alpha_{p}~;~\\\beta_{1}, ..., \beta_{q}~;~\end{array} z\right)
=\sum_{n=0}^{\infty}\frac{(\alpha_{1})_{n}  ...  (\alpha_{p})_{n}}{(\beta_{1})_{n} ... (\beta_{q})_{n}}\frac{z^{n}}{n!}~,\newline
\end{equation}
where $\alpha_{i}\in\mathbb{C}~(i=1,...,p)$ and $\beta_{j}\in\mathbb{C}\setminus \mathbb{Z}_{0}^{-}~(j=1,...,q)~\left(\ \mathbb{Z}_{0}^{-}:=\{0,-1,-2,...\}\right)$ and \\
$\left(\ p,~q\in\mathbb{N}_{0}:=\mathbb{N}\cup\{0\}=\{0,1,2,...\}\right)$.
The ${}_{p}F_{q}$ series in eq.(\ref{FAR4}) is convergent for $|z|<\infty$ if $p\leq q$, and for $|z|<1$ if $p=q+1$.
Furthermore, if we set
\begin{equation}\label{FAR6}
\omega=\left(\ \sum_{j=1}^{q}\beta_{j}-\sum_{i=1}^{p}\alpha_{i}\right),\newline
\end{equation}
it is known that the ${}_{p}F_{q}$ series, with $p=q+1$, is\\
(i) absolutely convergent for $|z|=1$ if $\Re(\omega)>0,$\\
(ii) conditionally convergent for $|z|=1, z\neq1$, if $-1<Re(\omega)\leq0$.\\
\\
If we replace $z$ by $\frac{z}{b}$ in the eq.(\ref{FAR3}) and taking the $\ell imit {|b|\rightarrow\infty}$ , we get  Kummer's confluent hypergeometric function, represented by\\
\begin{equation}\label{FAR7}
\ell imit_{|b|\rightarrow\infty}~{}_{2}F_{1}\left(\begin{array}{lll}a,~b~;\\~~~~ d~;\end{array} \frac{z}{b} \right)=
{}_{1}F_{1}\left(a~ ;~d;~z\right)=\sum_{n=0}^{\infty}\frac{(a)_{n}}{(d)_{n}}\frac{z^{n}}{n!},
\end{equation}
~~~~~where $|z|<\infty,~a\in\mathbb{C},~d\in\mathbb{C}\backslash \mathbb{Z}_{0}^{-}$.\\
When $d=b$ in the eq.(\ref{FAR3}) , we get a binomial function, given by
\begin{equation}\label{FAR8}
(1-z)^{-a}={}_{1}F_{0}\left(\begin{array}{lll}a~;\\  \overline{~~~};\end{array} z \right)
=\sum_{n=0}^{\infty}\frac{(a)_{n}}{n!}z^{n},
\end{equation}
~~~~~where $|z|<1,~~a\in\mathbb{C}$.\\
Next we collect some results that we will need in the sequel.\\
 Classical Gauss summation theorem for the series ${}_{2}F_{1}(\cdot)$ (\cite{B2},\cite{R}) is given by
\begin{equation}\label{FAR9}
{}_{2}F_{1}\left(\begin{array}{lll}a,~b~;\\~~~~ d~;\end{array} 1 \right)
=\frac{\Gamma(d)\Gamma(d-a-b)}{\Gamma(d-a)\Gamma(d-b)},
\end{equation}
~~~~~where $\Re(d-a-b)>0,~~d\in\mathbb{C}\backslash \mathbb{Z}_{0}^{-}$.\\
If $z=-1$ and $d=1+a-b$ in the eq.(\ref{FAR3}), we get Kummer's first summation theorem \cite{B2}, given by
\begin{equation}\label{FAR10}
{}_{2}F_{1}\left(\begin{array}{lll}a,~~b~~~~~~~;\\  1+a-b;\end{array} -1\right)
=\frac{\Gamma(1+a-b)\Gamma(1+\frac{a}{2})}{\Gamma(1+\frac{a}{2}-b)\Gamma(1+a)},
\end{equation}
~~~~~provided  $\Re(b)<1$,~~$1+a-b\in\mathbb{C}\backslash \mathbb{Z}_{0}^{-}$.\\
Also the classical summation theorem for hypergeometric series ${}_{4}F_{3}(-1)$ \cite[ p.28, Eq.(4.4.3)]{B2} is given by
\begin{equation}\label{FAR11}
{}_{4}F_{3}\left(\begin{array}{lll}a,~~1+\frac{a}{2},~~ b,~~c~~~~~~~~~~~~~;\\ \frac{a}{2},~~1+a-b,~~ 1+a-c;~\end{array} -1\right)
=\frac{\Gamma(1+a-b)\Gamma(1+a-c)}{\Gamma(1+a)\Gamma(1+a-b-c)}\newline,
\end{equation}
~~~~~provided $\Re(a-2b-2c)>-2 ;~ \frac{a}{2}, 1+a-b, 1+a-c\in\mathbb{C}\backslash \mathbb{Z}_{0}^{-}$.\\
\\
For every positive integer $m$ \cite[p.22, eq.(26)]{s1}, we have\\
\begin{equation}\label{FAR12}
  (\lambda)_{mn}=m^{mn}\prod_{j=1}^{m}\left(\frac{\lambda+j-1}{m}\right)_{n}~~~~~~~~;m\in\mathbb{N},~n\in\mathbb{N}_{0}\newline.
\end{equation}
When $m=2$, we get
\begin{equation}\label{FAR13}
(\lambda)_{2n}=2^{2n}\left(\frac{\lambda}{2}\right)_{n}\left(\frac{\lambda+1}{2}\right)_{n}\newline.
\end{equation}
From the above result (\ref{FAR12}), we get\\
\begin{equation}\label{FAR14}
  \Gamma(mz)=(2\pi)^{\frac{(1-m)}{2}}m^{mz-\frac{1}{2}}\prod_{j=1}^{m}\Gamma\left(z+\frac{j-1}{m}\right),~~~~mz\in\mathbb{C}\backslash \mathbb{Z}_{0}^{-},
\end{equation}
which is known as Gauss-Legendre multiplication theorem for Gamma function.
 When we put $m=2$ in the eq.(\ref{FAR14}), we get
\begin{equation}\label{FAR15}
\sqrt{\pi}\Gamma(2z)=2^{2z-1}\Gamma(z)\Gamma\left(z+\frac{1}{2}\right),~~~~~~~~~~2z\in\mathbb{C}\backslash \mathbb{Z}_{0}^{-},
\end{equation}
~~~~which is known as Legendre's duplication formula.\\
Algebraic property of Pochhammer symbol:
\begin{equation}\label{FAR16}
(\lambda)_{m+n}=(\lambda)_{m}(\lambda+m)_{n}=(\lambda)_{n}(\lambda+n)_{m}.
\end{equation}
Recurrence relation:
\begin{equation}\label{FAR17}
\Gamma(z+1)=z~\Gamma(z).
\end{equation}
Relation between circular function and Gamma function:
\begin{equation}\label{FAR18}
\Gamma(z)\Gamma(1-z)=\frac{\pi }{\sin (\pi z)}~;~z\neq0,\pm1,\pm2,\pm3,....
\end{equation}
If given function $f(t)$ is well defined for all real values of $t>0$,
then Laplace transform of $f(t)$, denoted by $g(p)$, and is given by the integral
\begin{equation}\label{FAR19}
\mathcal{L}[f(t);p]=\int_{0}^{\infty}e^{-pt}f(t)dt=g(p);~~~~\Re(p)>0,
\end{equation}
where $p$ is a complex variable and inverse Laplace transform of $g(p)$ is given by
\begin{equation}\label{FAR20}
\mathcal{L}^{-1}[g(p);t]=\frac{1}{2\pi i}\int_{c-i\infty}^{c+i\infty}e^{pt}g(p)dp=f(t),
\end{equation}
where $c$ is a real constant that exceeds the real part of all the singularities of $g(p)$.\\
Laplace transform of any constant $k$ is given by
\begin{equation}\label{FAR21}
\mathcal{L}[k;q]=\int_{0}^{\infty}e^{-qt}k~dt=\frac{k}{q},
\end{equation}
~~~~~provided
\begin{equation}\label{FAR22}
\Re(q)>0.
\end{equation}
Laplace transform of $t^{z-1}$ \cite[p.12, eq.(33)]{E1} is given by
\begin{equation}\label{FAR23}
\mathcal{L}[t^{z-1};S]=\int_{0}^{\infty}~e^{-St}t^{z-1}dt=\frac{\Gamma(z)}{S^{z}},
\end{equation}
~~~~where $Re(S)>0, 0<\Re(z)<\infty ~~OR~~ \Re(S)=0, 0<\Re(z)<1.$\\
The Laplace transforms of sine and cosine functions are recorded in the table \cite[p.150-154, Entry (1), Entry (43)]{E2} and are given by
\begin{equation}\label{FAR24}
\int_{0}^{\infty}~e^{-px}\sin(\alpha x)dx=\frac{\alpha}{\alpha^{2}+p^{2}},
\end{equation}
~~~~~~where $\Re(p)>|Im(\alpha)|$,\\
and
\begin{equation}\label{FAR25}
\int_{0}^{\infty}~e^{-px}\cos(\alpha x)dx=\frac{p}{\alpha^{2}+p^{2}},
\end{equation}
~~~~~~where $\Re(p)>|Im(\alpha)|$.\\
The Binomial expansion is given by
\begin{equation}\label{FAR26}
(A+B)^{N}=\sum_{r=0}^{N}\binom{N}{r}A^{N-r}~B^{r}=\sum_{r=0}^{N}\binom{N}{r}A^{r}~B^{N-r},
\end{equation}
where $N$ is positive integer and  Binomial coefficient holds the property:
\begin{equation}\label{FAR27}
\binom{N}{r}=\binom{N}{N-r}.
\end{equation}
The finite series representations of positive integral powers of hyperbolic functions:
 \begin{equation}\label{FAR28}
 \sum_{i=0}^{m-1}\bigg[\binom{2m}{i}\cosh \{(2m-2i)\beta x\}\bigg]+\frac{1}{2}\binom{2m}{m}=\frac{1}{2}\{e^{\beta x}+e^{-\beta x}\}^{2m}=2^{2m-1}\cosh^{2m}(\beta x),
 \end{equation}
 \begin{equation}\label{FAR29}
 \sum_{j=0}^{n-1}\bigg[(-1)^{j}\binom{2n}{j}\cosh \{(2n-2j)\gamma x\}\bigg]+\frac{(-1)^{n}}{2}\binom{2n}{n}=\frac{1}{2}\{e^{\gamma x}-e^{-\gamma x}\}^{2n}=2^{2n-1}\sinh^{2n}(\gamma x),
 \end{equation}
 where $m$ and $n$ are positive integer,
 \begin{equation}\label{FAR30}
 \sum_{k=0}^{p}\bigg[(-1)^{k}\binom{2p+1}{k}\sinh \{(2p+1-2k)\lambda x\}\bigg]=\frac{1}{2}\{e^{\lambda x}-e^{-\lambda x}\}^{2p+1}=2^{2p}\sinh^{2p+1}(\lambda x),
 \end{equation}
 \begin{equation}\label{FAR31}
 \sum_{\ell=0}^{q}\bigg[\binom{2q+1}{\ell}\cosh \{(2q+1-2\ell)\mu x\}\bigg]=\frac{1}{2}\{e^{\mu x}+e^{-\mu x}\}^{2q+1}=2^{2q}\cosh^{2q+1}(\mu x),
 \end{equation}
 ~~~~~~where $p$ and $q$ are non negative integer.\\
 \textbf{Proof}: We take right hand side of the eq.(\ref{FAR28}) and applying Binomial expansion (\ref{FAR26})
  \begin{multline}\label{FAR32}
\frac{1}{2}\{e^{\beta x}+e^{-\beta x}\}^{2m}=\frac{1}{2}\sum_{i=0}^{2m}\bigg[\binom{2m}{i}(e^{\beta x})^{2m-i}~(e^{-\beta x})^{i}\bigg]
=\frac{1}{2}\sum_{i=0}^{2m}\bigg[T_{i}\bigg],\\
=\frac{1}{2}\bigg[(T_0+T_{2m})+(T_1+T_{2m-1})+(T_2+T_{2m-2})+...+\\
+(T_{m-3}+T_{m+3})+(T_{m-2}+T_{m+2})+(T_{m-1}+T_{m+1})+T_{m}\bigg].
\end{multline}
Now putting the values of $T_0,T_1,T_2,...,T_{m-3},T_{m-2},T_{m-1},T_m,T_{m+1},T_{m+2},T_{m+3},...,T_{2m-2},T_{2m-1},T_{2m}$ in the above eq.(\ref{FAR32}). Applying again Binomial expansion (\ref{FAR26}) and Binomial property  (\ref{FAR27}), after simplifications we get the left hand side of the eq.(\ref{FAR28}). Similarly proof of eq.(\ref{FAR29}) is akin to that of eq.(\ref{FAR28}).
Also, taking right hand side of the eq.(\ref{FAR30}), given by
\begin{multline}\label{FAR33}
 \frac{1}{2}\{e^{\lambda x}-e^{-\lambda x}\}^{2p+1}=\frac{1}{2}\sum_{k=0}^{2p+1}\bigg[\binom{2p+1}{k}(e^{\lambda x})^{k}~(-e^{-\lambda x})^{2p+1-k}\bigg]
=\frac{1}{2}\sum_{k=0}^{2p+1}\bigg[U_{k}\bigg],\\
=\frac{1}{2}\bigg[(U_0+U_{2p+1})+(U_1+U_{2p})+(U_2+U_{2p-1})+...+\\
+(U_{p-2}+U_{p+3})+(U_{p-1}+U_{p+2})+(U_{p}+U_{p+1})\bigg].
 \end{multline}
 Now putting the values of
 $U_0,U_1,U_2,...,U_{p-2},U_{p-1},U_p,U_{p+1},U_{p+2},U_{p+3},...,U_{2p-1},U_{2p},U_{2p+1}$ in the above eq.(\ref{FAR33}). Applying again Binomial expansion (\ref{FAR26}) and Binomial property  (\ref{FAR27}), after simplifications we get the left hand side of  eq.(\ref{FAR30}). Similarly proof of eq.(\ref{FAR31}) is akin to that of the eq.(\ref{FAR30}).\\
 \\
 Relation between circular and hyperbolic functions:
 \begin{eqnarray}\label{FAR34}
  \sin(i\theta)= i~ \sinh\theta, ~~~~\sinh(i\theta)&=& i~ \sin\theta,
\end{eqnarray}
\begin{eqnarray}\label{FAR35}
  \cos(i\theta)=\cosh\theta, ~~~~ \cosh(i\theta)&=& \cos\theta.
 \end{eqnarray}
 The identities for the product of sine and cosine functions, are given by
 \begin{equation}\label{FAR36}
 2\sin(A)\sin(B)=\cos(A-B)-\cos(A+B),
 \end{equation}
 \begin{equation}
 2\cos(A)\cos(B)=\cos(A-B)+\cos(A+B),
 \end{equation}
 \begin{equation}\label{FAR37}
 2\sin(A)\cos(B)=\sin(A-B)+\sin(A+B).
 \end{equation}
 By the successive applications of the above formulas (\ref{FAR36})-(\ref{FAR37}), we can derive the following product formulas
 \begin{multline}\label{FAR001}
 \cos(A)\cos(B)\cos(C)\\=\frac{1}{4}\bigg[\cos(A-B-C)+\cos(A-B+C)+\cos(A+B-C)+\cos(A+B+C)\bigg],
 \end{multline}
 \begin{multline}
 \sin(A)\cos(B)\cos(C)\\=\frac{1}{4}\bigg[\sin(A-B-C)+\sin(A-B+C)+\sin(A+B-C)+\sin(A+B+C)\bigg],
 \end{multline}
 \begin{multline}
 \sin(A)\sin(B)\cos(C)\\=\frac{1}{4}\bigg[\cos(A-B-C)+\cos(A-B+C)-\cos(A+B-C)-\cos(A+B+C)\bigg],
 \end{multline}
   \begin{multline}
 \sin(A)\sin(B)\sin(C)\\=\frac{1}{4}\bigg[-\sin(A-B-C)+\sin(A-B+C)+\sin(A+B-C)-\sin(A+B+C)\bigg],
 \end{multline}
 \begin{multline}
 \cos(A)\cos(B)\cos(C)\cos(D)\\
 =\frac{1}{8}\bigg[\cos(A-B-C+D)+\cos(A-B-C-D)+\cos(A-B+C-D)+\cos(A-B+C+D)\\
 +\cos(A+B-C+D)+\cos(A+B-C-D)+\cos(A+B+C-D)+\cos(A+B+C+D)\bigg].
 \end{multline}
 \begin{multline}\label{FAR002}
 \sin(A)\cos(B)\cos(C)\cos(D)\\
 =\frac{1}{8}\bigg[\sin(A-B-C+D)+\sin(A-B-C-D)+\sin(A-B+C-D)+\sin(A-B+C+D)\\
 +\sin(A+B-C+D)+\sin(A+B-C-D)+\sin(A+B+C-D)+\sin(A+B+C+D)\bigg],
 \end{multline}
 \begin{multline}
 \sin(A)\sin(B)\cos(C)\cos(D)\\
 =\frac{1}{8}\bigg[\cos(A-B-C+D)+\cos(A-B-C-D)+\cos(A-B+C-D)+\cos(A-B+C+D)\\
 -\cos(A+B-C+D)-\cos(A+B-C-D)-\cos(A+B+C-D)-\cos(A+B+C+D)\bigg],
 \end{multline}
 \begin{multline}
 \sin(A)\sin(B)\sin(C)\cos(D)\\
 =\frac{1}{8}\bigg[-\sin(A-B-C+D)-\sin(A-B-C-D)+\sin(A-B+C-D)+\sin(A-B+C+D)\\
 +\sin(A+B-C+D)+\sin(A+B-C-D)-\sin(A+B+C-D)-\sin(A+B+C+D)\bigg],
 \end{multline}
  \begin{multline}\label{FAR00}
 \sin(A)\sin(B)\sin(C)\sin(D)\\
 =\frac{1}{8}\bigg[\cos(A-B-C+D)-\cos(A-B-C-D)+\cos(A-B+C-D)-\cos(A-B+C+D)\\
 -\cos(A+B-C+D)+\cos(A+B-C-D)-\cos(A+B+C-D)+\cos(A+B+C+D)\bigg],
 \end{multline}
$~~~$ In many papers of Moll et.al \cite{VM1,VM2,VM3, VM5, VM8, hyp, VM6,VMB,VMA,VMC,VM4,VM7, VMD}, which are published in the journal  \texttt{Scientia Series A: Mathematical Sciences}, where they evaluated many definite integrals of the table of Gradshteyn and Ryzhik, by using the change of independent variables. Our approach in this paper is to obtain the analytical solutions of some novel integrals, by using the hypergeometric approach, which is totally  different approach as given by  Moll et.al.\\
$~~~~~~$ The plan of this paper is as follows. First we obtain the Laplace transforms of an arbitrary power of some finite series containing hyperbolic sine and cosine functions in the \texttt{section 2} and the analytical solutions of some novel integrals are shown in \texttt{section 3}. Also, the Laplace transforms of positive integral powers of sine and cosine functions, and their
combinations of product are shown in \texttt{section 4}. Moreover, special cases of some novel integrals are given in \texttt{section 5}.

\section{Laplace transforms of an arbitrary power of some finite series containing hyperbolic sine and cosine functions having different arguments }
The following known results (\ref{FAR38})-(\ref{FAR41}) are given in the tables \cite[p.163, Entry (5), Entry (6)]{E2}; see also  \cite[pp.384-385, Entry(3.541)(1), Entry (3.542)(1)]{G}
and \cite[p.11, eq.(25)]{E1}
\begin{equation}\label{FAR38}
\int_{0}^{\infty}e^{-s x} [\cosh(\gamma x)-1]^{\nu}dx=\frac{2\nu}{2^{\nu}(s+\nu\gamma)}~B\left(\frac{s}{\gamma}-\nu,2\nu\right),
\end{equation}
\begin{equation}\label{FAR39}
~~~~~~~~~~~~~~~~~~~~~~~~~~~~~~~~~~~~~~=\frac{1}{2^{\nu}(\gamma)}~B\left(\frac{s}{\gamma}-\nu,2\nu+1\right),
\end{equation}
~~~~~where $\Re(\gamma)>0,~\Re(\nu)>-\frac{1}{2},~\Re(s)>\Re(\gamma\nu),~\frac{s}{\gamma}-\nu+1\in\mathbb{C}\backslash\mathbb{Z}_{0}^{-}$,
\begin{equation}\label{FAR40}
\int_{0}^{\infty}e^{-s x} [\sinh(\lambda x)]^{\nu}dx=\frac{\nu}{2^{\nu}(s+\lambda\nu)}~B\left(\frac{s}{2\lambda}-\frac{\nu}{2},\nu\right),
\end{equation}
\begin{equation}\label{FAR41}
~~~~~~~~~~~~~~~~~~~~~~~~~~~~~~~~~~~~~~~=\frac{1}{2^{1+\nu}(\lambda)}~B\left(\frac{s}{2\lambda}-\frac{\nu}{2},1+\nu\right),
\end{equation}
~~~~~where $\Re(\lambda)>0,~\Re(\nu)>-1,~\Re(s)>\Re(\lambda\nu),~\frac{s}{2\lambda}-\frac{\nu}{2}+1\in\mathbb{C}\backslash \mathbb{Z}_{0}^{-}$.\\
In this section, our work is motivated by the above known results (\ref{FAR38})-(\ref{FAR41}). The analytical solutions of the following Laplace transforms of arbitrary powers of some finite series, in terms of hypergeometric functions, holds true.\\
I. The arbitrary power of first finite series containing hyperbolic cosine function and its Laplace transform.
\begin{eqnarray}\label{FAR42}
\int_{0}^{\infty}e^{-s x}\bigg[\sum_{i=0}^{m-1}\binom{2m}{i}\cosh \{(2m-2i)\beta x\}+\frac{1}{2}\binom{2m}{m}\bigg]^{\nu}dx,~~~~~~~~~~~~~~~~~~~~~\nonumber\\
=\frac{1}{2^{\nu}(s-2m\nu\beta)}{}_{2}F_{1}\left(\begin{array}{lll}-2m\nu,~\frac{s}{2\beta}-m\nu ;\\  \frac{s}{2\beta}-m\nu+1~~~~~;\end{array} -1\right),
\end{eqnarray}
~~~~where $\Re(m\nu)>-1,~\Re(s-2m\beta\nu)>0,~\Re(\beta)>0,~\frac{s}{2\beta}-m\nu+1\in\mathbb{C}\backslash \mathbb{Z}_{0}^{-}$ and $m$ is positive integer.\\
II. The arbitrary power of second finite series containing hyperbolic cosine function and its Laplace transform.
\begin{eqnarray}\label{FAR43}
\int_{0}^{\infty}e^{-s x}\bigg[\sum_{j=0}^{n-1}(-1)^{j}\binom{2n}{j}\cosh \{(2n-2j)\gamma x\}+\frac{(-1)^{n}}{2}\binom{2n}{n}\bigg]^{\nu}dx,~~~~~~\nonumber\\
=\frac{1}{2^{\nu}(s-2n\nu\gamma)}{}_{2}F_{1}\left(\begin{array}{lll}-2n\nu,~\frac{s}{2\gamma}-n\nu ;\\  \frac{s}{2\gamma}-n\nu+1~~~~~;\end{array} 1\right),\nonumber\\
=\frac{2^{1-\nu}(n\nu)}{(s+2n\nu\gamma)}~B\left(\frac{s}{2\gamma}-n\nu,~2n\nu\right),~~~~~~~~~~~~~~~~
\end{eqnarray}
where $\Re(2n\nu)>-1,~\Re(s-2n\gamma\nu)>0,~\Re(\gamma)>0,~\frac{s}{2\gamma}-n\nu+1\in\mathbb{C}\backslash \mathbb{Z}_{0}^{-}$ and $n$ is positive integer.\\
III. The arbitrary power of third finite series containing hyperbolic sine function and its Laplace transform.
\begin{eqnarray}\label{FAR44}
\int_{0}^{\infty}e^{-s x}\bigg[\sum_{k=0}^{p}(-1)^{k}\binom{2p+1}{k}\sinh \{(2p+1-2k)\lambda x\}\bigg]^{\nu}dx,~~~~~~~~~~~~~~~~~~~~~\nonumber\\
=\frac{1}{2^{\nu}(s-2p\nu\lambda-\lambda\nu)}{}_{2}F_{1}\left(\begin{array}{lll}-2p\nu-\nu,~\frac{s}{2\lambda}-p\nu-\frac{\nu}{2} ;\\  \frac{s}{2\lambda}-p\nu-\frac{\nu}{2}+1~~~~~~~~~~~~;\end{array} 1\right),\nonumber\\
=\frac{(2p\nu+\nu)}{2^{\nu}(s+2p\nu\lambda+\nu\lambda)}~B\left(\frac{s}{2\lambda}-p\nu-\frac{\nu}{2},~2p\nu+\nu\right),~~~~~~~~~~~
\end{eqnarray}
where $\Re(2p\nu+\nu)>-1,~\Re(s-2p\lambda\nu-\lambda\nu)>0,~\Re(\lambda)>0,~\frac{s}{2\lambda}-p\nu-\frac{\nu}{2}+1\in\mathbb{C}\backslash \mathbb{Z}_{0}^{-}$ and $p$ is non-negative integer.\\
IV. The arbitrary power of fourth finite series containing hyperbolic cosine function and its Laplace transform.
\begin{eqnarray}\label{FAR45}
\int_{0}^{\infty}e^{-s x}\bigg[\sum_{\ell=0}^{q}\binom{2q+1}{\ell}\cosh \{(2q+1-2\ell)\mu x\}\bigg]^{\nu}dx,~~~~~~~~~~~~~~~~~~~~~\nonumber\\
=\frac{1}{2^{\nu}(s-2q\nu\mu-\mu \nu)}{}_{2}F_{1}\left(\begin{array}{lll}-2q\nu-\nu,~\frac{s}{2\mu}-q\nu-\frac{\nu}{2} ;\\  \frac{s}{2\mu}-q\nu-\frac{\nu}{2}+1~~~~~~~~~~~~;\end{array} -1\right),
\end{eqnarray}
where $\Re(2q\nu+\nu)>-2,~\Re(s-2q\mu\nu-\mu\nu)>0,~\Re(\mu)>0, \frac{s}{2\mu}-q\nu-\frac{\nu}{2}+1\in\mathbb{C}\backslash \mathbb{Z}_{0}^{-}$ and $q$ is non-negative integer.\\
\textbf{Proof}: Taking $\nu$th power on both sides of the equations (\ref{FAR28})-(\ref{FAR31}). Then  multiply both sides of the resulting equations by $e^{-sx}$ and integrate with respect to $x$ over the interval $(0,\infty)$. Then finally using Laplace transform formulas (\ref{FAR24}) and (\ref{FAR25}), we get the results stated in (\ref{FAR42})-(\ref{FAR45}).

\section{ Some novel integrals with suitable convergence conditions, in terms of Beta functions}
Many authors have studied some definite integrals containing the integrands as a quotient
of hyperbolic functions.  Mainly,  V. H. Moll et.al evaluated some definite integrals given in the table of Gradshteyn and Ryzhik, by using the change of independent variables. We obtain the analytical solutions of the some definite integrals, using hypergeometric approach.\\
\\
V. The fifth novel integral states that
\begin{eqnarray}\label{FAR46}
\int_{0}^{\infty}\frac{\cosh(2\alpha t)}{[\cosh(pt)]^{2\beta}}dt= 4^{\beta-1}p^{-1} B\left(\beta+\frac{\alpha}{p},\beta-\frac{\alpha}{p}\right),
\end{eqnarray}
~~~~~where $\Re(\beta)<1, ~\Re(p)>0,~\Re\left(\beta\pm\frac{\alpha}{p}\right)>0$ and $\beta, \beta\pm\frac{\alpha}{p}+1\in\mathbb{C}\backslash \mathbb{Z}_{0}^{-}$.\\
\\
VI. The sixth novel integral states that
\begin{equation}\label{FAR47}
\int_{0}^{\infty}\frac{[\sinh(x)]^{\alpha}}{[\cosh(x)]^{\beta}}dx=\frac{1}{2}B\left(\frac{1+\alpha}{2},\frac{\beta-\alpha}{2}\right),
\end{equation}
~~~~~where $\Re(\alpha)>-1,~ -2<\Re(\alpha-\beta)<0$ and $\frac{\beta\pm\alpha+2}{2}\in\mathbb{C}\backslash \mathbb{Z}_{0}^{-}$.\\
\\
VII. The seventh novel integral states that
\begin{eqnarray}\label{FAR48}
\int_{0}^{\infty}\frac{\cos(ax)}{[\cosh(\beta x)]^{\nu}}dx
=\frac{2^{\nu-2}}{\beta\Gamma(\nu)}\Gamma\left(\frac{\nu}{2}+\frac{ia}{2\beta}\right)\Gamma\left(\frac{\nu}{2}-\frac{ia}{2\beta}\right),
\end{eqnarray}
\begin{equation}\label{FAR49}
=B\left(\frac{\nu}{2}+\frac{ia}{2\beta},\frac{\nu}{2}-\frac{ia}{2\beta}\right),
\end{equation}
~~~~~where $\Re(\beta)>0,~\Re(\nu)<2,~\Re(\nu\beta\pm ia)>0$ and $\frac{\nu}{2},\frac{\nu}{2}\pm\frac{ia}{2\beta}+1\in\mathbb{C}\backslash \mathbb{Z}_{0}^{-}.$\\
\\
\textbf{Hypergeometric proof of integral (\ref{FAR46})}: Suppose left hand side of eq.(\ref{FAR46}) is denoted by $\mho(\alpha,\beta,p)$. Using the well known result of hyperbolic function (\ref{FAR1}) and binomial function (\ref{FAR8}), in L.H.S. of  the above eq.(\ref{FAR46}), we obtain\\
\begin{equation*}\label{FAR50}
\mho(\alpha,\beta,p)=2^{2\beta-1}\int_{0}^{\infty}\bigg[e^{-2\beta pt}(e^{2\alpha t}+e^{-2\alpha t})
{}_{1}F_{0}\left(\begin{array}{lll}2\beta;\\ \overline{~~~};\end{array} -e^{-2p t}\right)\bigg]dt,
\end{equation*}
when $\Re(p)>0$ , then $|-e^{-2pt}|<1$ for all $t>0$. It is the convergence conditions of above Binomial function  ${}_{1}F_{0}(\cdot)$.
\begin{equation}\label{FAR51}
\mho(\alpha,\beta,p)=2^{2\beta-1}\int_{0}^{\infty}\bigg[e^{-2\beta pt}(e^{2\alpha t}+e^{-2\alpha t})\sum_{r=0}^{\infty}\frac{(2\beta)_{r}(-1)^{r}}{r!}e^{-2prt}\bigg]dt.
\end{equation}
Change the order of integration and the summation in eq.(\ref{FAR51}), we get
\begin{equation}\label{FAR52}
\mho(\alpha,\beta,p)=2^{2\beta-1}\sum_{r=0}^{\infty}\bigg[\frac{(2\beta)_{r}(-1)^{r}}{r!}\int_{0}^{\infty}\left(e^{-(2\beta p-2\alpha+2pr)t}+e^{-(2\beta p+2\alpha+2pr)t}\right)dt\bigg],
\end{equation}
where $\Re\left(\beta\pm\frac{\alpha}{p}\right)>0$, it is the convergence conditions of Laplace transform of unity in the integral (\ref{FAR52}). Then
using Laplace formula (\ref{FAR21}) in the above eq.(\ref{FAR52}). which yields
\begin{equation*}\label{FAR53}
\mho(\alpha,\beta,p)=2^{2\beta-1}\sum_{r=0}^{\infty}\bigg[\frac{(2\beta)_{r}(-1)^{r}}{r!}\left\{ \frac{p\beta+pr}{(\beta p+pr)^{2}-\alpha^{2}}\right\}\bigg],
\end{equation*}
\begin{equation}\label{FAR54}
=2^{2\beta-1}\sum_{r=0}^{\infty}\bigg[\frac{(2\beta)_{r}(-1)^{r}}{r!}\frac{\{\beta p+(p)r\}}{[(\beta p+\alpha)+pr][(\beta p-\alpha)+pr]}\bigg].
\end{equation}
Using algebraic properties of Pochhammer symbol, after simplifications, we obtain
\begin{equation*}\label{FAR55}
\mho(\alpha,\beta,p)=2^{2\beta-1}\left(\frac{\beta p}{\beta^{2}p^{2}-\alpha^{2}}\right)\sum_{r=0}^{\infty}\bigg[\frac{(2\beta)_{r}(\beta+1)_{r}(\beta+\frac{\alpha}{p})_{r}(\beta-\frac{\alpha}{p})_{r}(-1)^{r}}{r!(\beta)_{r}(\beta+\frac{\alpha}{p}+1)_{r}(\beta-\frac{\alpha}{p}+1)_{r}}\bigg],
\end{equation*}
\begin{equation}\label{FAR56}
=2^{2\beta-1}\left(\frac{\beta p}{\beta^{2}p^{2}-\alpha^{2}}\right)~~{}_{4}F_{3}\left(\begin{array}{lll} 2\beta,\beta+1, \beta+\frac{\alpha}{p}, \beta-\frac{\alpha}{p};\\\beta, \beta+\frac{\alpha}{p}+1, \beta-\frac{\alpha}{p}+1;\end{array} -1\right),
\end{equation}
~~where $\Re(\beta)<1$  and  $\beta, \beta\pm\frac{\alpha}{p}+1\in\mathbb{C}\backslash \mathbb{Z}_{0}^{-}$; it is the convergence conditions of  summation  ${}_{4}F_{3}(-1)$ in the eq.(\ref{FAR56}). Then
applying the classical summation theorem (\ref{FAR11}), in  the eq.(\ref{FAR56}), which yields
\begin{equation}\label{FAR57}
\mho(\alpha,\beta,p)=2^{2\beta-1}\left(\frac{\beta p}{\beta^{2}p^{2}-\alpha^{2}}\right)\frac{\Gamma(\beta+\frac{\alpha}{p}+1)\Gamma(\beta-\frac{\alpha}{p}+1)}{\Gamma(1+2\beta)}.
\end{equation}
In view of the recurrence relation (\ref{FAR17}), we get
\begin{equation}\label{FAR58}
\mho(\alpha,\beta,p)=4^{\beta-1}p^{-1} ~\frac{\Gamma(\beta+\frac{\alpha}{p})\Gamma(\beta-\frac{\alpha}{p})}{\Gamma(2\beta)}\newline.
\end{equation}
Finally, using the Beta function in terms of Gamma function (\ref{FAR2}) in the above eq. (\ref{FAR58}), we get R.H.S. of the result (\ref{FAR46}).\\
\textbf{Hypergeometric proof of integral (\ref{FAR47})}:  Suppose left hand side of eq.(\ref{FAR47}) is denoted by $\Lambda(\alpha,\beta)$. Applying the result (\ref{FAR1}) in the left hand side of the eq.(\ref{FAR47}), we get
\begin{equation*}\label{FAR59}
\Lambda(\alpha,\beta)=2^{(\beta-\alpha)}\int_{0}^{\infty}\bigg[e^{(\alpha-\beta)x}{}_{1}F_{0}\left(\begin{array}{lll}-\alpha;\\ \overline{~~~~};\end{array} e^{-2x}\right)
{}_{1}F_{0}\left(\begin{array}{lll}\beta;\\ \overline{~~~};\end{array} -e^{-2x}\right)\bigg]dx,
\end{equation*}
where $|\pm e^{-2x}|<1$ for all $x>0$,
\begin{equation}\label{FAR60}
\Lambda(\alpha,\beta)=2^{(\beta-\alpha)}\int_{0}^{\infty}\bigg[e^{(\alpha-\beta)x}\sum_{r=0}^{\infty}\sum_{s=0}^{\infty}\frac{(-\alpha)_{r}(\beta)_{s}(e^{-2x})^{r}(e^{-2x})^{s}(-1)^{s}}{r!s!}\bigg]dx.
\end{equation}
Change the order of integration and the summation in eq.(\ref{FAR60}), which yields
\begin{equation}\label{FAR61}
\Lambda(\alpha,\beta)=2^{(\beta-\alpha)}\sum_{r=0}^{\infty}\sum_{s=0}^{\infty}\bigg[\frac{(-\alpha)_{r}(\beta)_{s}(-1)^{s}}{r!s!}
\int_{0}^{\infty}e^{-\{(\beta-\alpha)+2(r+s)\}x}dx\bigg],
\end{equation}
where $\Re(\beta-\alpha)>0$, it is the convergence conditions of Laplace transform of unity in the integral (\ref{FAR61}).Then
applying formula (\ref{FAR21}) in the eq.(\ref{FAR61}), we obtain
\begin{equation}\label{FAR62}
\Lambda(\alpha,\beta)=2^{(\beta-\alpha)}\sum_{r=0}^{\infty}\sum_{s=0}^{\infty}\bigg[\frac{(-\alpha)_{r}(\beta)_{s}(-1)^{s}}{s!~r!}\left\{\frac{1}{(\beta-\alpha)+2(r+s)}\right\}\bigg].
\end{equation}
Using algebraic properties of Pochhammer symbol in the eq.(\ref{FAR62}), after  simplifications, we get
\begin{equation*}\label{FAR63}
\Lambda(\alpha,\beta)=\frac{2^{(\beta-\alpha)}}{(\beta-\alpha)}\sum_{r=0}^{\infty}\sum_{s=0}^{\infty}\bigg[\frac{(-\alpha)_{r}(\beta)_{s}\left(\frac{\beta-\alpha}{2}\right)_{s}\left(\frac{\beta-\alpha}{2}+s\right)_{r}(-1)^{s}}
{\left(\frac{\beta-\alpha+2}{2}\right)_{s}\left(\frac{\beta-\alpha+2}{2}+s\right)_{r}~s!~r!}\bigg],
\end{equation*}
\begin{equation}\label{FAR64}
=\frac{2^{(\beta-\alpha)}}{(\beta-\alpha)}\sum_{s=0}^{\infty}\bigg[\frac{(\beta)_{s}\left(\frac{\beta-\alpha}{2}\right)_{s}(-1)^{s}}{\left(\frac{\beta-\alpha+2}{2}\right)_{s}~s!}
{}_{2}F_{1}\left(\begin{array}{lll}-\alpha,~\frac{\beta-\alpha}{2}+s;~\\  \frac{\beta-\alpha+2}{2}+s;\end{array} 1\right)\bigg],
\end{equation}
where $\frac{\beta-\alpha+2}{2}\in\mathbb{C}\backslash \mathbb{Z}_{0}^{-}$,~~$\Re(1+\alpha)>0$; it is the convergence conditions of above classical Gauss summation theorem  ${}_{2}F_{1}(1)$ in the eq.(\ref{FAR64}). Then
using the result (\ref{FAR9}) in the eq.(\ref{FAR64}), we get
\begin{equation*}\label{FAR65}
\Lambda(\alpha,\beta)=\frac{2^{(\beta-\alpha)}}{(\beta-\alpha)}\frac{\Gamma(1+\alpha)\Gamma\left(\frac{\beta-\alpha+2}{2}\right)}{\Gamma\left(\frac{\beta+\alpha+2}{2}\right)}
\sum_{s=0}^{\infty}\bigg[\frac{(\beta)_{s}\left(\frac{\beta-\alpha}{2}\right)_{s}(-1)^{s}}{\left(\frac{\beta+\alpha+2}{2}\right)_{s}~s!}\bigg],
\end{equation*}
\begin{equation}\label{FAR66}
=\frac{2^{(\beta-\alpha)}}{(\beta-\alpha)}\frac{\Gamma(1+\alpha)\Gamma(\frac{\beta-\alpha+2}{2})}{\Gamma(\frac{\beta+\alpha+2}{2})}
{}_{2}F_{1}\left(\begin{array}{lll}\beta, \frac{\beta-\alpha}{2};\\  \frac{\beta+\alpha+2}{2};\end{array} -1\right),
\end{equation}
~where $\frac{\beta+\alpha+2}{2}\in\mathbb{C}\backslash \mathbb{Z}_{0}^{-}$,~ $\Re(\alpha-\beta)>-2$; it is the convergence conditions of Kummer's first summation theorem  ${}_{2}F_{1}(-1)$ in the eq.(\ref{FAR66}). Then using the result (\ref{FAR10}) in the eq.(\ref{FAR66}), we obtain
\begin{equation}\label{FAR67}
\Lambda(\alpha,\beta)=\frac{2^{(\beta-\alpha)}}{(\beta-\alpha)}\frac{\sqrt{\pi}~\Gamma(1+\alpha)\Gamma\left(\frac{\beta-\alpha}{2}+1\right)\Gamma\left(1+\frac{\beta}{2}\right)}{\sqrt{\pi}\Gamma(1+\beta)\Gamma\left(1+\frac{\alpha}{2}\right)}.
\end{equation}
 Applying  Legendre's duplication formula (\ref{FAR15}) in the eq. (\ref{FAR67}), we get
 \begin{equation}\label{FAR68}
 \Lambda(\alpha,\beta)=\frac{\Gamma\left(\frac{1+\alpha}{2}\right)\Gamma\left(\frac{\beta-\alpha}{2}\right)}{2\Gamma\left(\frac{1+\beta}{2}\right)}.
 \end{equation}
 Lastly, using Beta function in terms of Gamma function (\ref{FAR2}) in the above eq. (\ref{FAR68}), we get right hand side of the result (\ref{FAR47}).\\
 \textbf{Hypergeometric proof of integral (\ref{FAR48})}: Suppose left hand side of eq.(\ref{FAR48}) is denoted by $\Phi(a,\beta,\nu)$. Using the well known result of hyperbolic function (\ref{FAR1}) and binomial function (\ref{FAR8}), in L.H.S. of  the above eq.(\ref{FAR48}), we obtain\\
\begin{equation*}\label{FAR69}
\Phi(a,\beta,\nu)=2^{\nu-1}\int_{0}^{\infty}\bigg[e^{-\nu\beta x}(e^{iax}+e^{-iax})
{}_{1}F_{0}\left(\begin{array}{lll}\nu;\\ \overline{~~~};\end{array} -e^{-2\beta x}\right)\bigg]dx,
\end{equation*}
when $\Re(\beta)>0$, then $|-e^{-2\beta x}|<1$ for all $x>0$,
\begin{equation}\label{FAR70}
\Phi(a,\beta,\nu)=2^{\nu-1}\int_{0}^{\infty}\bigg[e^{-\nu\beta x}(e^{iax}+e^{-iax})\sum_{r=0}^{\infty}\frac{(\nu)_{r}(-1)^{r}}{r!}e^{-2\beta rx}\bigg]dx.
\end{equation}
Change the order of integration and the summation, we get
\begin{equation}\label{FAR71}
\Phi(a,\beta,\nu)=2^{\nu-1}\sum_{r=0}^{\infty}\bigg[\frac{(\nu)_{r}(-1)^{r}}{r!}\int_{0}^{\infty}\left(e^{-(\nu\beta-ia+2\beta r)x}+e^{-(\nu\beta +ia+2\beta r)x}\right)dx\bigg],
\end{equation}
where $\Re(\nu\beta\pm ia)>0$, it is the convergence conditions of Laplace transform of unity in the eq.(\ref{FAR71}).
Then applying formula (\ref{FAR21}) in the eq.(\ref{FAR71}) and using algebraic properties of Pochhammer symbol, after  simplifications, we obtain
\begin{equation*}\label{FAR72}
\Phi(a,\beta,\nu)=\frac{2^{\nu}(\nu\beta) }{(\beta^{2}\nu^{2}+a^{2})}
\sum_{r=0}^{\infty}\bigg[\frac{(\nu)_{r}(\frac{\nu}{2}+1)_{r}(\frac{\nu}{2}-\frac{ia}{2\beta})_{r}(\frac{\nu}{2}+\frac{ia}{2\beta})_{r}(-1)^{r}}{r!(\frac{\nu}{2})_{r}(\frac{\nu}{2}-\frac{ia}{2\beta}+1)_{r}(\frac{\nu}{2}+\frac{ia}{2\beta}+1)_{r}}\bigg],
\end{equation*}
\begin{equation}\label{FAR73}
=\frac{2^{\nu}(\nu\beta) }{(\beta^{2}\nu^{2}+a^{2})}~~{}_{4}F_{3}\left(\begin{array}{lll} \nu,\frac{\nu}{2}+1, \frac{\nu}{2}-\frac{ia}{2\beta}, \frac{\nu}{2}+\frac{ia}{2\beta};\\ \frac{\nu}{2}, \frac{\nu}{2}-\frac{ia}{2\beta}+1, \frac{\nu}{2}+\frac{ia}{2\beta}+1;\end{array} -1\right),
\end{equation}
where $\Re(\nu)<2$ and $\frac{\nu}{2},\frac{\nu}{2}\pm\frac{ia}{2\beta}+1\in\mathbb{C}\backslash \mathbb{Z}_{0}^{-} $; it is the convergence conditions of the function ${}_{4}F_{3}(-1)$ in the eq.(\ref{FAR73}). Then using summation theorem ${}_{4}F_{3}(-1)$ in the eq.(\ref{FAR73}), we get right hand side of the result (\ref{FAR48}).
\section{Laplace transforms of positive integral powers of sine and cosine functions having different arguments, and their combinations of product}
 The following results (\ref{FAR74})-(\ref{FAR81}) are motivated by the work given in the tables \cite[pp.150-156, Entry(3),Entry(7),Entry(47),Entry(51)]{G3}; see also  \cite[sec.(3.611-4.146)]{G}.
Each of the following Laplace transforms of positive integral powers of sine and cosine functions having different arguments, holds true:
\begin{equation}\label{FAR74}
\int_{0}^{\infty}~e^{-sx}\cos^{2m}(\beta x)dx=
\frac{1}{2^{2m-1}}\sum_{i=0}^{m-1}\bigg[\binom{2m}{i}\frac{s}{\{(2m-2i)\beta\}^{2}+s^{2}}\bigg]+\frac{1}{2^{2m}}\binom{2m}{m}\frac{1}{s},
\end{equation}
\begin{equation}\label{FAR75}
~~~~~~~~~=\frac{1}{2^{2m}\{s-\omega~2m\beta\}}~~{}_{2}F_{1}\left(\begin{array}{lll} -2m,~\frac{-\omega s-2m\beta}{2\beta}~;\\\frac{-\omega s-2m\beta}{2\beta}+1~~~~~;\end{array} -1\right),
\end{equation}
where $\Re(s)>2m|Im(\beta)|$, $\omega=\sqrt{(-1)}$ and $m$ is positive integer,
\begin{equation}\label{FAR76}
\int_{0}^{\infty}~e^{-sx}\sin^{2n}(\gamma x)dx=
\frac{(-1)^{n}}{2^{2n-1}}~\sum_{j=0}^{n-1}\bigg[(-1)^{j}\binom{2n}{j}\frac{s}{\{(2n-2j)\gamma\}^{2}+s^{2}}\bigg]+\frac{1}{2^{2n}}\binom{2n}{n}\frac{1}{s},
\end{equation}
\begin{equation}\label{FAR77}
  =\frac{(2n)!}{2^{2n}\{s+2\omega n\gamma\}\left(1+\frac{\omega s}{2\gamma}\right)_{n}\left(-\frac{\omega s}{2\gamma}\right)_{n}},~~~~~
\end{equation}
where $\Re(s)>2n|Im(\gamma)|$ , $\omega=\sqrt{(-1)}$ and $n$ is positive integer,
\begin{equation}\label{FAR78}
\int_{0}^{\infty}~e^{-sx}\sin^{2p+1}(\lambda x)dx=
\frac{(-1)^{p}}{2^{2p}}\sum_{k=0}^{p}\bigg[(-1)^{k}\binom{2p+1}{k}\frac{(2p+1-2k)\lambda}{\{(2p+1-2k)\lambda\}^{2}+s^{2}}\bigg],
\end{equation}
\begin{equation}\label{FAR79}
  =\frac{(2p+1)!~\lambda}{2^{2p}\{s-\omega\lambda\}\{s+\omega\lambda(2p+1)\}\left(\frac{3\lambda+\omega s }{2\lambda}\right)_{p}\left(\frac{\lambda-\omega s}{2\lambda}\right)_{p}},~~~~~
\end{equation}
where $\Re(s)>(2p+1)|Im(\lambda)|$ , $\omega=\sqrt{(-1)}$ and $p$ is non negative integer, \\
\begin{equation}\label{FAR80}
\int_{0}^{\infty}~e^{-sx}\cos^{2q+1}(\mu x)dx=
\frac{1}{2^{2q}}\sum_{\ell=0}^{q}\bigg[\binom{2q+1}{\ell}\frac{s}{\{(2q+1-2\ell)\mu\}^{2}+s^{2}}\bigg],
\end{equation}
\begin{equation}\label{FAR81}
~~~~~~~~~=\frac{1}{2^{2q+1}\{s-\omega~(2q+1)\mu\}}~~{}_{2}F_{1}\left(\begin{array}{lll} -2q-1,~\frac{-\omega s-(2q+1)\mu}{2\mu}~;\\\frac{-\omega s-(2q+1)\mu}{2\mu}+1~~~~~~~~~;\end{array} -1\right),
\end{equation}
where $\Re(s)>(2q+1)|Im(\mu)|$ and $\omega=\sqrt{(-1)}$ and $q$ is non negative integer.\\
\\
\textbf{Proof}: Multiply both sides of the following equations  (\ref{FAR82})-(\ref{FAR85}) by $e^{-sx}$ and integrate with respect to $x$ over the interval $(0,\infty)$. Then finally using Laplace transform formulas (\ref{FAR21}), (\ref{FAR24}) and (\ref{FAR25}), we get the results stated in (\ref{FAR74})-(\ref{FAR81})
 \begin{equation}\label{FAR82}
 \sum_{i=0}^{m-1}\bigg[\binom{2m}{i}\cos \{(2m-2i)\beta x\}\bigg]+\frac{1}{2}\binom{2m}{m}=2^{2m-1}\cos^{2m}(\beta x),
 \end{equation}
 \begin{equation}\label{FAR83}
 \sum_{j=0}^{n-1}\bigg[(-1)^{j}\binom{2n}{j}\cos \{(2n-2j)\gamma x\}\bigg]+\frac{(-1)^{n}}{2}\binom{2n}{n}=(-1)^{n}~2^{2n-1}\sin^{2n}(\gamma x),
 \end{equation}
 \begin{equation}\label{FAR84}
 \sum_{k=0}^{p}\bigg[(-1)^{k}\binom{2p+1}{k}\sin \{(2p+1-2k)\lambda x\}\bigg]=(-1)^{p}~2^{2p}\sin^{2p+1}(\lambda x),
 \end{equation}
 \begin{equation}\label{FAR85}
 \sum_{\ell=0}^{q}\bigg[\binom{2q+1}{\ell}\cos \{(2q+1-2\ell)\mu x\}\bigg]=2^{2q}\cos^{2q+1}(\mu x).
 \end{equation}
 Using the relations  (\ref{FAR34}) and (\ref{FAR35}) in the equations (\ref{FAR28})-(\ref{FAR31}) with suitable adjustment of parameters, we can obtain the above results (\ref{FAR82})-(\ref{FAR85}).
The Laplace transforms of the products of two functions, three functions and four function  having different arguments and different powers taken together, obtained with the help of above four equations (\ref{FAR82})-(\ref{FAR85}), can be written by\\
\\
$\bullet$ Ten integrals associated with  Laplace transforms of products of two functions:
\begin{equation}\label{FAR86}
\int_{0}^{\infty}~e^{-sx}\sin^{2m}(\beta x)\cos^{2n}(\gamma x)dx,
\end{equation}
\begin{equation}
\int_{0}^{\infty}~e^{-sx}\sin^{2m+1}(\beta x)\cos^{2n+1}(\gamma x)dx,
\end{equation}
\begin{equation}
\int_{0}^{\infty}~e^{-sx}\sin^{2m}(\beta x)\cos^{2n+1}(\gamma x)dx,
\end{equation}
\begin{equation}
\int_{0}^{\infty}~e^{-sx}\sin^{2m+1}(\beta x)\cos^{2n}(\gamma x)dx,
\end{equation}
\begin{equation}
\int_{0}^{\infty}~e^{-sx}\sin^{2m}(\beta x)\sin^{2n}(\gamma x)dx,
\end{equation}
\begin{equation}
\int_{0}^{\infty}~e^{-sx}\sin^{2m+1}(\beta x)\sin^{2n+1}(\gamma x)dx,
\end{equation}
\begin{equation}
\int_{0}^{\infty}~e^{-sx}\sin^{2m}(\beta x)\sin^{2n+1}(\gamma x)dx,
\end{equation}
\begin{equation}
\int_{0}^{\infty}~e^{-sx}\cos^{2m}(\beta x)\cos^{2n}(\gamma x)dx,
\end{equation}
\begin{equation}
\int_{0}^{\infty}~e^{-sx}\cos^{2m+1}(\beta x)\cos^{2n+1}(\gamma x)dx,
\end{equation}
\begin{equation}
\int_{0}^{\infty}~e^{-sx}\cos^{2m}(\beta x)\cos^{2n+1}(\gamma x)dx,
\end{equation}
$\bullet$ Twenty integrals associated with Laplace transforms of products of three functions:
\begin{equation}
\int_{0}^{\infty}~e^{-sx}\sin^{2m}(\beta x)\cos^{2n}(\gamma x)\sin^{2p}(\lambda x)dx,
\end{equation}
\begin{equation}
\int_{0}^{\infty}~e^{-sx}\sin^{2m}(\beta x)\cos^{2n}(\gamma x)\sin^{2p+1}(\lambda x)dx,
\end{equation}
\begin{equation}
\int_{0}^{\infty}~e^{-sx}\sin^{2m}(\beta x)\cos^{2n}(\gamma x)\cos^{2p}(\lambda x)dx,
\end{equation}
\begin{equation}
\int_{0}^{\infty}~e^{-sx}\sin^{2m}(\beta x)\cos^{2n}(\gamma x)\cos^{2p+1}(\lambda x)dx,
\end{equation}
\begin{equation}
\int_{0}^{\infty}~e^{-sx}\sin^{2m+1}(\beta x)\cos^{2n+1}(\gamma x)\sin^{2p}(\lambda x)dx,
\end{equation}
\begin{equation}
\int_{0}^{\infty}~e^{-sx}\sin^{2m+1}(\beta x)\cos^{2n+1}(\gamma x)\sin^{2p+1}(\lambda x)dx,
\end{equation}
\begin{equation}
\int_{0}^{\infty}~e^{-sx}\sin^{2m+1}(\beta x)\cos^{2n+1}(\gamma x)\cos^{2p}(\lambda x)dx,
\end{equation}
\begin{equation}
\int_{0}^{\infty}~e^{-sx}\sin^{2m+1}(\beta x)\cos^{2n+1}(\gamma x)\cos^{2p+1}(\lambda x)dx,
\end{equation}
\begin{equation}
\int_{0}^{\infty}~e^{-sx}\sin^{2m}(\beta x)\cos^{2p+1}(\gamma x)\sin^{2n}(\lambda x)dx,
\end{equation}
\begin{equation}
\int_{0}^{\infty}~e^{-sx}\sin^{2m}(\beta x)\cos^{2n+1}(\gamma x)\cos^{2p+1}(\lambda x)dx,
\end{equation}
\begin{equation}
\int_{0}^{\infty}~e^{-sx}\sin^{2m+1}(\beta x)\cos^{2n}(\gamma x)\sin^{2p+1}(\lambda x)dx,
\end{equation}
\begin{equation}
\int_{0}^{\infty}~e^{-sx}\sin^{2m+1}(\beta x)\cos^{2n}(\gamma x)\cos^{2p}(\lambda x)dx,
\end{equation}
\begin{equation}
\int_{0}^{\infty}~e^{-sx}\sin^{2m}(\beta x)\sin^{2n}(\gamma x)\sin^{2p}(\lambda x)dx,
\end{equation}
\begin{equation}
\int_{0}^{\infty}~e^{-sx}\sin^{2m}(\beta x)\sin^{2n}(\gamma x)\sin^{2p+1}(\lambda x)dx,
\end{equation}
\begin{equation}
\int_{0}^{\infty}~e^{-sx}\sin^{2m+1}(\beta x)\sin^{2n+1}(\gamma x)\sin^{2p}(\lambda x)dx,
\end{equation}
\begin{equation}
\int_{0}^{\infty}~e^{-sx}\sin^{2m+1}(\beta x)\sin^{2n+1}(\gamma x)\sin^{2p+1}(\lambda x)dx,
\end{equation}
\begin{equation}
\int_{0}^{\infty}~e^{-sx}\cos^{2m}(\beta x)\cos^{2n}(\gamma x)\cos^{2p}(\lambda x)dx,
\end{equation}
\begin{equation}
\int_{0}^{\infty}~e^{-sx}\cos^{2m}(\beta x)\cos^{2n}(\gamma x)\cos^{2p+1}(\lambda x)dx,
\end{equation}
\begin{equation}
\int_{0}^{\infty}~e^{-sx}\cos^{2m+1}(\beta x)\cos^{2n+1}(\gamma x)\cos^{2p}(\lambda x)dx,
\end{equation}
\begin{equation}
\int_{0}^{\infty}~e^{-sx}\cos^{2m+1}(\beta x)\cos^{2n+1}(\gamma x)\cos^{2p+1}(\lambda x)dx,
\end{equation}

$\bullet$ Thirty six integrals associated with Laplace transforms of products of four functions:
\begin{equation}
\int_{0}^{\infty}~e^{-sx}\sin^{2m}(\beta x)\cos^{2n}(\gamma x)\sin^{2p}(\lambda x)\sin^{2q}(\mu x)dx,
\end{equation}
\begin{equation}
\int_{0}^{\infty}~e^{-sx}\sin^{2m}(\beta x)\cos^{2n}(\gamma x)\sin^{2p}(\lambda x)\sin^{2q+1}(\mu x)dx,
\end{equation}
\begin{equation}
\int_{0}^{\infty}~e^{-sx}\sin^{2m}(\beta x)\cos^{2n}(\gamma x)\sin^{2p}(\lambda x)\cos^{2q}(\mu x)dx,
\end{equation}
\begin{equation}
\int_{0}^{\infty}~e^{-sx}\sin^{2m}(\beta x)\cos^{2n}(\gamma x)\sin^{2p}(\lambda x)\cos^{2q+1}(\mu x)dx,
\end{equation}
\begin{equation}
\int_{0}^{\infty}~e^{-sx}\sin^{2m}(\beta x)\cos^{2n}(\gamma x)\sin^{2p+1}(\lambda x)\sin^{2q+1}(\mu x)dx,
\end{equation}
\begin{equation}
\int_{0}^{\infty}~e^{-sx}\sin^{2m}(\beta x)\cos^{2n}(\gamma x)\sin^{2p+1}(\lambda x)\cos^{2q}(\mu x)dx,
\end{equation}
\begin{equation}
\int_{0}^{\infty}~e^{-sx}\sin^{2m}(\beta x)\cos^{2n}(\gamma x)\sin^{2p+1}(\lambda x)\cos^{2q+1}(\mu x)dx,
\end{equation}
\begin{equation}
\int_{0}^{\infty}~e^{-sx}\sin^{2m}(\beta x)\cos^{2n}(\gamma x)\cos^{2p}(\lambda x)\cos^{2q}(\mu x)dx,
\end{equation}
\begin{equation}
\int_{0}^{\infty}~e^{-sx}\sin^{2m}(\beta x)\cos^{2n}(\gamma x)\cos^{2p}(\lambda x)\cos^{2q+1}(\mu x)dx,
\end{equation}
\begin{equation}
\int_{0}^{\infty}~e^{-sx}\sin^{2m}(\beta x)\cos^{2n}(\gamma x)\cos^{2p+1}(\lambda x)\cos^{2q+1}(\mu x)dx,
\end{equation}
\begin{equation}
\int_{0}^{\infty}~e^{-sx}\sin^{2m+1}(\beta x)\cos^{2n+1}(\gamma x)\sin^{2p}(\lambda x)\sin^{2q}(\mu x)dx,
\end{equation}
\begin{equation}
\int_{0}^{\infty}~e^{-sx}\sin^{2m+1}(\beta x)\cos^{2n+1}(\gamma x)\sin^{2p}(\lambda x)\sin^{2q+1}(\mu x)dx,
\end{equation}
\begin{equation}
\int_{0}^{\infty}~e^{-sx}\sin^{2m+1}(\beta x)\cos^{2n+1}(\gamma x)\sin^{2p}(\lambda x)\cos^{2q+1}(\mu x)dx,
\end{equation}
\begin{equation}
\int_{0}^{\infty}~e^{-sx}\sin^{2m+1}(\beta x)\cos^{2n+1}(\gamma x)\sin^{2p+1}(\lambda x)\sin^{2q+1}(\mu x)dx,
\end{equation}
\begin{equation}
\int_{0}^{\infty}~e^{-sx}\sin^{2m+1}(\beta x)\cos^{2n+1}(\gamma x)\sin^{2p+1}(\lambda x)\cos^{2q}(\mu x)dx,
\end{equation}
\begin{equation}
\int_{0}^{\infty}~e^{-sx}\sin^{2m+1}(\beta x)\cos^{2n+1}(\gamma x)\sin^{2p+1}(\lambda x)\cos^{2q+1}(\mu x)dx,
\end{equation}
\begin{equation}
\int_{0}^{\infty}~e^{-sx}\sin^{2m+1}(\beta x)\cos^{2n+1}(\gamma x)\cos^{2p}(\lambda x)\cos^{2q}(\mu x)dx,
\end{equation}
\begin{equation}
\int_{0}^{\infty}~e^{-sx}\sin^{2m+1}(\beta x)\cos^{2n+1}(\gamma x)\cos^{2p}(\lambda x)\cos^{2q+1}(\mu x)dx,
\end{equation}
\begin{equation}
\int_{0}^{\infty}~e^{-sx}\sin^{2m+1}(\beta x)\cos^{2n+1}(\gamma x)\cos^{2p+1}(\lambda x)\cos^{2q+1}(\mu x)dx,
\end{equation}
\begin{equation}
\int_{0}^{\infty}~e^{-sx}\sin^{2m}(\beta x)\cos^{2n+1}(\gamma x)\sin^{2p}(\lambda x)\sin^{2q}(\mu x)dx,
\end{equation}
\begin{equation}
\int_{0}^{\infty}~e^{-sx}\sin^{2m}(\beta x)\cos^{2n+1}(\gamma x)\sin^{2p}(\lambda x)\sin^{2q+1}(\mu x)dx,
\end{equation}
\begin{equation}
\int_{0}^{\infty}~e^{-sx}\sin^{2m}(\beta x)\cos^{2n+1}(\gamma x)\sin^{2p}(\lambda x)\cos^{2q+1}(\mu x)dx,
\end{equation}
\begin{equation}
\int_{0}^{\infty}~e^{-sx}\sin^{2m}(\beta x)\cos^{2n+1}(\gamma x)\cos^{2p+1}(\lambda x)\cos^{2q+1}(\mu x)dx,
\end{equation}
\begin{equation}
\int_{0}^{\infty}~e^{-sx}\sin^{2m+1}(\beta x)\cos^{2n}(\gamma x)\sin^{2p+1}(\lambda x)\sin^{2q+1}(\mu x)dx,
\end{equation}
\begin{equation}
\int_{0}^{\infty}~e^{-sx}\sin^{2m+1}(\beta x)\cos^{2n}(\gamma x)\sin^{2p+1}(\lambda x)\cos^{2q}(\mu x)dx,
\end{equation}
\begin{equation}
\int_{0}^{\infty}~e^{-sx}\sin^{2m+1}(\beta x)\cos^{2n}(\gamma x)\cos^{2p}(\lambda x)\cos^{2q}(\mu x)dx,
\end{equation}
\begin{equation}
\int_{0}^{\infty}~e^{-sx}\sin^{2m}(\beta x)\sin^{2n}(\gamma x)\sin^{2p}(\lambda x)\sin^{2q}(\mu x)dx,
\end{equation}
\begin{equation}
\int_{0}^{\infty}~e^{-sx}\sin^{2m}(\beta x)\sin^{2n}(\gamma x)\sin^{2p}(\lambda x)\sin^{2q+1}(\mu x)dx,
\end{equation}
\begin{equation}
\int_{0}^{\infty}~e^{-sx}\sin^{2m}(\beta x)\sin^{2n}(\gamma x)\sin^{2p+1}(\lambda x)\sin^{2q+1}(\mu x)dx,
\end{equation}
\begin{equation}
\int_{0}^{\infty}~e^{-sx}\sin^{2m+1}(\beta x)\sin^{2n+1}(\gamma x)\sin^{2p}(\lambda x)\sin^{2q+1}(\mu x)dx,
\end{equation}
\begin{equation}
\int_{0}^{\infty}~e^{-sx}\sin^{2m+1}(\beta x)\sin^{2n+1}(\gamma x)\sin^{2p+1}(\lambda x)\sin^{2q+1}(\mu x)dx,
\end{equation}
\begin{equation}
\int_{0}^{\infty}~e^{-sx}\cos^{2m}(\beta x)\cos^{2n}(\gamma x)\cos^{2p}(\lambda x)\cos^{2q}(\mu x)dx,
\end{equation}
\begin{equation}
\int_{0}^{\infty}~e^{-sx}\cos^{2m}(\beta x)\cos^{2n}(\gamma x)\cos^{2p}(\lambda x)\cos^{2q+1}(\mu x)dx,
\end{equation}
\begin{equation}
\int_{0}^{\infty}~e^{-sx}\cos^{2m}(\beta x)\cos^{2n}(\gamma x)\cos^{2p+1}(\lambda x)\cos^{2q+1}(\mu x)dx,
\end{equation}
\begin{equation}
\int_{0}^{\infty}~e^{-sx}\cos^{2m+1}(\beta x)\cos^{2n+1}(\gamma x)\cos^{2p}(\lambda x)\cos^{2q+1}(\mu x)dx,
\end{equation}
\begin{equation}\label{FAR87}
\int_{0}^{\infty}~e^{-sx}\cos^{2m+1}(\beta x)\cos^{2n+1}(\gamma x)\cos^{2p+1}(\lambda x)\cos^{2q+1}(\mu x)dx,
\end{equation}
The above combinations of Laplace transforms (\ref{FAR86})-(\ref{FAR87}) (not recorded in the tables \cite[Vol.1, sec.(4.7), pp.150-160]{E2}; \cite{G3,P5}) can be derived by putting finite series representations of the corresponding even and odd positive integral powers of sine and cosine functions with the help of equations  (\ref{FAR82})-(\ref{FAR85}). In the resulting product of finite summations, use identities (\ref{FAR36})-(\ref{FAR00}), then apply  the Laplace transformation formulas (\ref{FAR21}), (\ref{FAR24}) and  (\ref{FAR25}).
\section{Special cases of novel integrals of sections 2 and 3, with suitable convergence conditions}
The following integrals from (\ref{FAR97})-(\ref{FAR106}) are not found in the literature of integral transforms\\
$\bullet$ If we set $m=1$ and $\beta$ is replaced by $\beta/2$ in the eq.(\ref{FAR42}), we get
\begin{equation}\label{FAR97}
\int_{0}^{\infty}e^{-s x} [\cosh(\beta x)+1]^{\nu}dx=
\frac{1}{2^{\nu}(s-\beta\nu)}{}_{2}F_{1}\left(\begin{array}{lll}-2\nu,~\frac{s}{\beta}-\nu;\\  \frac{s}{\beta}-\nu+1;\end{array} -1\right),
\end{equation}
~~~~~where $\Re(\beta)>0,~\Re(\nu)>-1,~\Re(s)>\Re(\beta\nu)$,~$ \frac{s}{\beta}-\nu+1\in\mathbb{C}\backslash \mathbb{Z}_{0}^{-}$.\\
\\
$\bullet$ If we set $m=2$ and $\beta$ is replaced by $\beta/4$ in the eq.(\ref{FAR42}), we have
\begin{eqnarray}\label{FAR98}
\int_{0}^{\infty}e^{-s x} \left[\cosh(\beta x)+4\cosh\left(\frac{\beta x}{2}\right)+3\right]^{\nu}dx\nonumber~~~~~~~~~~~~~~~~~~~~~~~~~~~~~~\\
=\frac{1}{2^{\nu}(s-\nu\beta)}{}_{2}F_{1}\left(\begin{array}{lll}-4\nu,~\frac{2s}{\beta}-2\nu ;\\  \frac{2s}{\beta}-2\nu+1;\end{array} -1\right),
\end{eqnarray}
~~~~~where $\Re(\beta)>0,\Re(\nu)>-\frac{1}{2},~\Re(s)>\Re(\beta\nu),~\frac{2s}{\beta}-2\nu+1\in\mathbb{C}\backslash \mathbb{Z}_{0}^{-}$.\\
\\
$\bullet$ If we set $m=3$ and $\beta$ is replaced by $\beta/6$ in the eq.(\ref{FAR42}), we obtain
\begin{eqnarray}\label{FAR99}
\int_{0}^{\infty}e^{-s x} \left[\cosh(\beta x)+6\cosh\left(\frac{2\beta x}{3}\right)+15\cosh\left(\frac{\beta x}{3}\right)+10\right]^{\nu}dx,\nonumber\\
=\frac{1}{2^{\nu}(s-\nu\beta)}{}_{2}F_{1}\left(\begin{array}{lll}-6\nu,~\frac{3s}{\beta}-3\nu ;\\  \frac{3s}{\beta}-3\nu+1;\end{array} -1\right),
\end{eqnarray}
~~~~~where $\Re(\beta)>0,\Re(\nu)>-\frac{1}{3} ,~\Re(s)>\Re(\beta\nu),~\frac{3s}{\beta}-3\nu+1\in\mathbb{C}\backslash \mathbb{Z}_{0}^{-}$.\\
\\
$\bullet$ If we set $n=1$ and $\gamma$ is replaced by $\gamma/2$ in the eq.(\ref{FAR43}), we get the eq.(\ref{FAR38}).\\
\\
$\bullet$ If we set $n=2$ and $\gamma$ is replaced by $\gamma/4$ in the eq.(\ref{FAR43}), we obtain
\begin{equation}\label{FAR100}
\int_{0}^{\infty}e^{-s x} \left[\cosh(\gamma x)-4\cosh\left(\frac{\gamma x}{2}\right)+3\right]^{\nu}dx=\frac{4\nu}{2^{v}(s+\gamma\nu)}~B\left(\frac{2s}{\gamma}-2\nu,~4\nu\right),
\end{equation}
~~~~~where $\Re(\gamma)>0,~\Re(\nu)>-\frac{1}{4},~\Re(s)>\Re(\gamma\nu),~\frac{2s}{\gamma}-2\nu+1\in\mathbb{C}\backslash\mathbb{Z}_{0}^{-}$.\\
\\
$\bullet$ If we set $n=3$ and $\gamma$ is replaced by $\gamma/6$ in the eq.(\ref{FAR43}), we have
\begin{eqnarray}\label{FAR101}
\int_{0}^{\infty}e^{-s x} \left[\cosh(\gamma x)-6\cosh\left(\frac{2\gamma x}{3}\right)+15\cosh\left(\frac{\gamma x}{3}\right)-10\right]^{\nu}dx,\nonumber\\
=\frac{6\nu}{2^{v}(s+\gamma\nu)}~B\left(\frac{3s}{\gamma}-3\nu,~6\nu\right),~~~~~~~~~~~
\end{eqnarray}
~~~~~where $\Re(\gamma)>0,~\Re(\nu)>-\frac{1}{6},~\Re(s)>\Re(\gamma\nu),~\frac{3s}{\gamma}-3\nu+1\in\mathbb{C}\backslash \mathbb{Z}_{0}^{-}$.\\
\\
$\bullet$ If we set $p=0$ in the eq.(\ref{FAR44}), we get the eq.(\ref{FAR40}).
\\
\\
$\bullet$ If we set $p=1$ and $\lambda$ is replaced by $\lambda/3$ in the eq.(\ref{FAR44}), we have
\begin{eqnarray}\label{FAR102}
\int_{0}^{\infty}e^{-s x} \left[\sinh(\lambda x)-3\sinh\left(\frac{\lambda x}{3}\right)\right]^{\nu}dx,~~~~~~~~~~~~~~~~~~~~~~~~\nonumber\\
=\frac{3\nu}{2^{v}(s+\lambda\nu)}~B\left(\frac{3s}{2\lambda}-\frac{3\nu}{2},~3\nu\right),~~~~~~~~~~~
\end{eqnarray}
~~~~~where $\Re(\lambda)>0,~\Re(\nu)>-\frac{1}{3},~\Re(s)>\Re(\lambda\nu),~\frac{3s}{2\lambda}-\frac{3\nu}{2}+1\in\mathbb{C}\backslash \mathbb{Z}_{0}^{-}$.\\
\\
$\bullet$ If we set $p=2$ and $\lambda$ is replaced by $\lambda/5$ in the eq.(\ref{FAR44}), we obtain
\begin{eqnarray}\label{FAR103}
\int_{0}^{\infty}e^{-s x} \left[\sinh(\lambda x)-5\sinh\left(\frac{3\lambda x}{5}\right)+10\sinh\left(\frac{\lambda x}{5}\right)\right]^{\nu}dx,\nonumber\\
=\frac{5\nu}{2^{v}(s+\lambda\nu)}~B\left(\frac{5s}{2\lambda}-\frac{5\nu}{2},~5\nu\right),~~~~~~~~~~~
\end{eqnarray}
~~~~~where $\Re(\lambda)>0,~\Re(\nu)>-\frac{1}{5},~\Re(s)>\Re(\lambda\nu),~\frac{5s}{2\lambda}-\frac{5\nu}{2}+1\in\mathbb{C}\backslash \mathbb{Z}_{0}^{-}$.\\
\\
$\bullet$ If we set $q=0$ in the eq.(\ref{FAR45}), we get
\begin{equation}\label{FAR104}
\int_{0}^{\infty}e^{-s x} [\cosh(\mu x)]^{\nu}dx=\frac{1}{2^{\nu}(s-\mu\nu)}
{}_{2}F_{1}\left(\begin{array}{lll}-\nu,~\frac{s}{2\mu}-\frac{\nu}{2};\\  \frac{s}{2\mu}-\frac{\nu}{2}+1;\end{array} -1\right),
\end{equation}
~~~~~where $\Re(\mu)>0,~\Re(\nu)>-2,~\Re(s)>\Re(\mu\nu)$,~$ \frac{s}{2\mu}-\frac{\nu}{2}+1\in\mathbb{C}\backslash \mathbb{Z}_{0}^{-}$.\\
\\
$\bullet$ If we set $q=1$ and $\mu$ is replaced by $\mu/3$ in the eq.(\ref{FAR45}), we have
\begin{eqnarray}\label{FAR105}
\int_{0}^{\infty}e^{-s x} \left[\cosh(\mu x)+3\cosh\left(\frac{\mu x}{3}\right)\right]^{\nu}dx,~~~~~~~~~~~~~~~~~~~~~~~~~~\nonumber\\
=\frac{1}{2^{\nu}(s-\mu\nu)}
{}_{2}F_{1}\left(\begin{array}{lll}-3\nu,~\frac{3s}{2\mu}-\frac{3\nu}{2};\\  \frac{3s}{2\mu}-\frac{3\nu}{2}+1;\end{array} -1\right),
\end{eqnarray}
~~~~~where $\Re(\mu)>0,~\Re(\nu)>-\frac{2}{3},~\Re(s)>\Re(\mu\nu),~ \frac{3s}{2\mu}-\frac{3\nu}{2}+1\in\mathbb{C}\backslash \mathbb{Z}_{0}^{-}$.\\
\\
$\bullet$ If we set $q=2$ and $\mu$ is replaced by $\mu/5$ in the eq.(\ref{FAR45}), we obtain
\begin{eqnarray}\label{FAR106}
\int_{0}^{\infty}e^{-s x} \left[\cosh(\mu x)+5\cosh\left(\frac{3\mu x}{5}\right)+10\cosh\left(\frac{\mu x}{5}\right)\right]^{\nu}dx,\nonumber
\\
=\frac{1}{2^{\nu}(s-\mu\nu)}
{}_{2}F_{1}\left(\begin{array}{lll}-5\nu,~\frac{5s}{2\mu}-\frac{5\nu}{2};\\  \frac{5s}{2\mu}-\frac{5\nu}{2}+1;\end{array} -1\right),
\end{eqnarray}
~~~~~where $\Re(\mu)>0,~\Re(\nu)>-\frac{2}{5},~\Re(s)>\Re(\mu\nu),~ \frac{5s}{2\mu}-\frac{5\nu}{2}+1\in\mathbb{C}\backslash \mathbb{Z}_{0}^{-}$.\\
\\
 The following integrals from (\ref{FAR107})-(\ref{FAR116}) containing the integrands as a quotient of hyperbolic and trigonometric functions, are solved by Cauchy's residue theorem and are available in any standard text book of complex analysis or functions of a complex variable.\\
 \\
$\bullet$ If we put $p=1,~\alpha=0$ and $\beta=\mu$ in the above eq.(\ref{FAR46}), which yields
\begin{equation}\label{FAR107}
\int_{0}^{\infty}\frac{1}{\{\cosh(t)\}^{2\mu}}dt= 4^{\mu-1}B\left(\mu,\mu\right),~~~~\Re(\mu)>0.
\end{equation}
 $\bullet$ If we set $2\alpha=a,~p=b,~2\beta=1$ in the above eq.(\ref{FAR46}), we obtain
 \begin{equation}\label{FAR108}
 \int_{0}^{\infty}\frac{\cosh(at)}{\cosh(b t)}dt=\frac{\pi}{2b}\sec\left(\frac{a\pi}{2b}\right),~~~~~~~~~~~~~~b>|a|.
 \end{equation}
 $\bullet$ If we set $2\alpha=a,~p=\pi,~2\beta=1$ in the above eq.(\ref{FAR46}), we get
 \begin{equation}\label{FAR109}
 \int_{0}^{\infty}\frac{\cosh(at)}{\cosh(\pi t)}dt=\frac{1}{2}\sec\left(\frac{a}{2}\right),~~~~~~~~~~~~~~-\pi<a< \pi.
 \end{equation}
 $\bullet$ If we set $2\alpha=a,~p=1,~2\beta=1$ in the above eq.(\ref{FAR46}), we have
 \begin{equation}\label{FAR110}
 \int_{0}^{\infty}\frac{\cosh(at)}{\cosh( t)}dt=\frac{\pi}{2}\sec\left(\frac{a\pi}{2}\right),~~~~~~~~~~~~~~|a|<1.
 \end{equation}
  $\bullet$ If we set $\alpha=0,~2\beta=1$ in the above eq.(\ref{FAR46}), we obtain
  \begin{equation}\label{FAR111}
  \int_{0}^{\infty}\frac{dt}{\cosh(at)}=\frac{\pi}{2a},~~~~~~~~~~~~~~~~~~~~~~~~~~~~~a>0.
  \end{equation}
  $\bullet$ If we set $\beta=1$ in the above eq.(\ref{FAR46}), we have
   \begin{equation}\label{FAR112}
 \int_{0}^{\infty}\frac{\cosh(2\alpha t)}{\cosh^{2}(pt)}dt=\frac{\pi}{p}cosec\left(\frac{\pi \alpha}{p}\right).
   \end{equation}
 $\bullet$ If we set $\alpha=0$, $\beta=\mu$ in the eq. (\ref{FAR47}), we get
\begin{equation}\label{FAR113}
\int_{0}^{\infty}\frac{1}{\{\cosh(t)\}^{\mu}}dt=\frac{1}{2}B\left(\frac{1}{2},\frac{\mu}{2}\right),~~~~\Re(\mu)>0.
\end{equation}
 $\bullet$ If we set $\nu=2n$ in the above eq.(\ref{FAR48}), we obtain
\begin{eqnarray}
\int_{0}^{\infty}\frac{\cos(ax)}{[\cosh(\beta x)]^{2n}}dx
=\frac{2^{2n-2}}{\beta\Gamma(2n)}\Gamma\left(n+\frac{ia}{2\beta}\right)\Gamma\left(n-\frac{ia}{2\beta}\right),
\end{eqnarray}
~~~~~where $\Re(\beta)>0$ and $n$ is a positive integer.\\
\\
 $\bullet$ If we set $a=0,~\beta=1,~\nu=2$ in the above eq.(\ref{FAR48}), we get
  \begin{equation}\label{FAR114}
  \int_{0}^{\infty}\frac{1}{\cosh^{2}(x)}dx=1,
  \end{equation}
  $\bullet$ If we set $\nu=2$ in the above eq.(\ref{FAR48}), we have
  \begin{equation}\label{FAR115}
  \int_{0}^{\infty}\frac{\cos(ax)}{\cosh^{2}(\beta x)}dx=\frac{\pi a}{2\beta^{2}\sinh(\frac{\pi a}{2\beta})},
  \end{equation}
  ~~~where $\Re(\beta)>0,~a>0$.\\
  $\bullet$ If we set $\beta=1$ and differentiate with respect to $a$ in the above eq.(\ref{FAR115}). Then applying the Leibnitz's rule for differentiation under the sign of integration, we obtain
  \begin{equation}\label{FAR116}
  \int_{0}^{\infty}\frac{x\sin(ax)}{\cosh^{2}(x)}dx=\frac{2\pi\sinh(\frac{a\pi}{2})-a\pi^{2}\cosh(\frac{a\pi}{2})}{4\sinh^{2}(\frac{\pi a}{2})},
  \end{equation}
  ~~~where $~a>0$.
\section*{Conclusion}
Here, we have described some definite integrals containing the quotients of hyperbolic functions in terms of Beta functions. Also, Laplace transforms of an arbitrary power of some finite series containing hyperbolic sine and cosine functions and Laplace transforms of positive integral powers of sine and cosine functions,  in terms of Gauss hypergeometric  function and  Beta function. Thus certain integrals of hyperbolic and trigonometric functions, which may be different from those of presented here, can also be evaluated by the  hypergeometric approach.

\end{document}